\documentclass[10pt]{amsart}
\usepackage[english]{babel}
\usepackage{latexsym}

\DeclareMathOperator{\cent}{\bf Z}

\newcommand{\F}{\mathbb F}

\newcommand{\N}{\mathbb N}

\newtheorem{dummy}{Dummy}

\numberwithin{dummy}{section}
\numberwithin{equation}{section}

\newtheorem{lemma}[dummy]{Lemma}

\newtheorem{theorem}[dummy]{Theorem}

\newtheorem{prop}[dummy]{Proposition}

\theoremstyle{definition}

\theoremstyle{remark}
\newtheorem{rem}[dummy]{Remark}
\begin{document}%___________________________________________________________

\bibliographystyle{amsalpha}
\author{Marina Avitabile}
\author{Sandro Mattarei}
\thanks{Both authors are members of INdAM-GNSAGA, Italy.
They acknowledge financial  support from Ministero dell'Istruzione, dell'Universit\`a  e
  della  Ricerca, Italy, to  the
  project ``Graded Lie algebras  and pro-$p$-groups of finite width''.}
\title[Diamonds of finite type in thin {Lie} algebras]{Diamonds of finite type in thin {L}ie algebras}
\date{\today}
\begin{abstract}
Borrowing some terminology from pro-$p$ groups,
thin Lie algebras are $\N$-graded Lie algebras of width two and obliquity zero,
generated in degree one.
In particular, their homogeneous components have degree one or two, and they are termed {\em diamonds}
in the latter case.
In one of the two main subclasses of thin Lie algebras
the earliest diamond after that in degree one occurs in degree $2q-1$,
where $q$ is a power of the characteristic.
This paper is a contribution to a classification project of this subclass of thin Lie algebras.
Specifically, we prove that, under certain technical assumptions,
the degree of the earliest diamond {\em of finite type} in such a Lie algebra
can only have a certain form, which does occur in explicit examples constructed elsewhere.
\end{abstract}
%(CONTROLLARE!!!)}
\maketitle

\section{Introduction}\label{sec:introduction}%____________________________________________________

The coclass conjectures formulated in~\cite{L-GN} have had a considerable impact
on the theory of $p$-groups and pro-$p$ groups.
After the completion of their proofs in~\cite{L-G} and~\cite{Sha:coclass}
a study of a coclass theory for Lie algebras was initiated in~\cite{ShZe:finite-coclass}.
It is natural at a first stage to focus on Lie algebras over a field, which are graded
over the positive integers, and which are generated by their homogeneous component of degree one.
The results of~\cite{ShZe:finite-coclass} show that analogues of the coclass conjectures
hold for these Lie algebras, provided the underlying field has characteristic zero.
For our present purposes, the main assertion of the coclass conjectures (now theorems)
is that these Lie algebras are soluble,
and this implies some sort of periodicity in their structure (in the infinite-dimensional case).

Examples constructed in~\cite{Sha:max} show that in positive characteristic
analogues of the coclass conjectures fail even for the special case of coclass one, namely,
graded Lie algebras of maximal class need not be soluble.
(In the rest of this discussion we tacitly assume that graded Lie algebras are generated by their
homogeneous component of degree one, but results without this hypothesis were obtained in~\cite{ShZe:finite-coclass}
in characteristic zero and in~\cite{CV-L} in positive characteristic.)
The examples in~\cite{Sha:max} are (the positive parts of twisted) {\em loop algebras} of certain non-classical
finite-dimensional simple Lie algebras.
By their very construction, they do have a periodic structure, despite not being soluble.
This is not always the case for graded Lie algebras of maximal class.
In fact, the results of~\cite{CMN} show that the majority of graded Lie algebras of maximal class are
not periodic.
Subsequent work has lead to a classification of graded Lie algebras of maximal class
in arbitrary odd~\cite{CN} and even characteristic~\cite{Ju:maximal}.

The possibility of a classification of $p$-groups and pro-$p$ groups according to invariants
other than the coclass was suggested in~\cite{KL-GP}; see also Chapter~12 of~\cite{L-GMcKay}.
The three invariant suggested are called {\em width,} {\em obliquity} and {\em rank}.
Finite $p$-groups of width two and obliquity zero, which constitute the simplest non-trivial case,
had already been introduced in~\cite{Br,BCS} under the name of {\em thin} groups.
Analogous {\em thin} Lie algebras have been studied in several papers.

According to~\cite{CMNS}, a {\em thin} Lie algebra is a graded Lie algebra
$L= \bigoplus_{i=1}^{\infty} L_i$
over a field $\F$,
with $\dim(L_1)=2$ and satisfying the following
\emph{covering property}:
\begin{equation}\label{eq:CP}
\text{$L_{i+1}=[u,L_1]$ for every $0 \neq u \in L_i$, for all $i\ge 1$.}
\end{equation}
It was shown in~\cite{CMNS} that the covering property is equivalent to the following
analogue for Lie algebras of having {\em obliquity zero}:
every graded ideal $I$ of $L$ is located between two consecutive terms $L^i$ of the lower central series of $L$,
in the sense that $L^i\supseteq I\supseteq L^{i+1}$ for some $i$.
The definition of thin Lie algebra implies at once that $L$ is generated by $L_1$ as a Lie algebra,
and that every homogeneous component has dimension $1$ or $2$ (whence $L$ has {\em width} two).
We call a {\em diamond} any homogeneous component of dimension $2$,
and a {\em chain} the sequence of one-dimensional  homogeneous components between two consecutive diamonds.
In particular $L_1$ is a diamond, the {\em first diamond}.
If that is the only diamond, $L$ has maximal class.
However, it is convenient to explicitly exclude graded Lie algebras of maximal class from the definition
of thin Lie algebras, and hence we assume that a thin Lie algebra has at least two diamonds.

It is known from~\cite{CMNS} and~\cite{AviJur} that the second diamond of an infinite-dimensional thin Lie algebra
of zero or odd characteristic can only occur in degree $3$, $5$, $q$ or $2q-1$
where $q$ is a power of the characteristic $p$,
in case this is positive.
All these possibilities really do occur, and examples can be constructed by taking
{\em loop algebras} of suitable finite-dimensional Lie algebras.
In particular, examples of thin Lie algebras with second diamond in degree $3$ or $5$ can be produced
as loop algebras of classical Lie algebras of type $A_1$ or $A_2$ (see~\cite{CMNS} and~\cite{Car:thin_addendum})
and arise also as the graded Lie algebras associated to the lower central series of
certain $p$-adic analytic pro-$p$ groups of the corresponding types (see~\cite{Mat:thin-groups}, and also~\cite{KL-GP}).
It was also proved in~\cite{CMNS} that those loop algebras (two isomorphism types in case of $A_1$)
are the only infinite-dimensional thin Lie algebras
with second diamond in degree $3$ or $5$ (for $p>3$, and respectively $p>5$,
and with a further assumption in the former case).
In particular, although solubility fails, periodicity in the structure holds
for thin Lie algebras with second diamond in degree $3$ or $5$
(under the assumptions mentioned).

The remaining cases are typically modular, and most known examples
involve finite-dimensional simple Lie algebras of Cartan type.
In this paper, including the present discussion, we assume the characteristic of the ground field to be odd.
Thin Lie algebras with the second diamond in degree $q$ have been called {\em of Nottingham type}
in~\cite{Car:Nottingham,Car:Zassenhaus-three,Car:thin_addendum,CaMa:Nottingham,AviMat:Nottingham},
because the simplest example is the graded Lie algebra associated with the lower central series of the
so-called {\em Nottingham group}.
We refer to the Introduction of~\cite{AviMat:Nottingham} for an up-to-date overview of our knowledge
of this class of thin Lie algebras.

The present paper is concerned with thin Lie algebras having second diamond in degree $2q-1$,
which have been called $(-1)$-algebras in~\cite{CaMa:thin}.
Under certain assumptions, which we make precise in Section~\ref{sec:preliminaries},
each diamond of a $(-1)$-algebra can be assigned a {\em type} taking
value in the underlying field plus infinity.
Diamonds of type zero are really one-dimensional components,
but it is convenient to allow them in certain positions and dub them {\em fake};
we will add the attribute {\em genuine} to explicitly exclude that a diamond may be fake.
The locations and types of the diamonds constitute the most crucial
piece of information required to describe a $(-1)$-algebra.
In fact, it is conceivable that this piece of information may determine the algebra completely.
We discuss this further at the end of Section~\ref{sec:preliminaries},
and prove a partial result in Section~\ref{sec:chains} which supports this possibility,
namely, Proposition~\ref{prop:chains}.

Several diamond patterns are possible.
The $(-1)$-algebras with all diamonds of infinite type need not be periodic,
but they are well understood.
In fact, an explicit bijection was constructed in~\cite{CaMa:thin}
between them and a certain subclass of the graded Lie algebras of maximal class.

The $(-1)$-algebras with all diamonds of finite type were also constructed in~\cite{CaMa:thin}
and were shown to be determined by a certain finite-dimensional quotient.
More precisely, there is a unique infinite-dimensional thin Lie algebra
with second diamond in degree $2q-1$ and of finite type:
its diamonds occur at regular intervals (provide we include some fake diamonds in the count)
and have all finite types which follow an arithmetic progression.
An explicit construction for this thin Lie algebras as a loop algebra of a certain
finite-dimensional simple Lie algebras is given in~\cite{CaMa:Hamiltonian}.

Finally, there are $(-1)$-algebras with diamonds of both finite and infinite types.
According to the previous paragraph, the second diamond here must be of infinite type,
and a diamond of finite type must occur later.
These $(-1)$-algebras are constructed in~\cite{AviMat:-1}, and have diamonds in all degrees
congruent to $1$ modulo $q-1$, except those congruent to $q$ modulo $p^{s+1}$
(for some $s\ge 1$).
If we are willing to accept as diamonds of type zero (that is, fake)
the one-dimensional components in degrees congruent to $q$ modulo $p^{s+1}(q-1)$,
then there are diamonds in all degrees congruent to $1$ modulo $q-1$.
The diamonds occur in sequences of $p^s-1$ diamonds of infinite type
separated by single occurrences of diamonds of finite type,
and the latter types follow an arithmetic progression.

We believe that the Lie algebras described in the previous three paragraphs
should exhaust the class of $(-1)$-algebras (at least in odd characteristic).
While the results of~\cite{CaMa:thin} mentioned above take care of the $(-1)$-algebras
with all diamonds of infinite type, extensive machine calculations suggest
that the occurrence of a single diamond of finite type in a $(-1)$-algebra
should force a sort of periodicity and, therefore, cause the algebra to be uniquely determined.
According to~\cite{CaMa:thin} this is the case if the second diamond
has finite type (and hence, as it turns out, all diamonds have finite type).
This uniqueness result in~\cite{CaMa:thin} was proved by
showing that the homogeneous relations of degree up to $2q$ satisfied by any thin Lie
algebra with second diamond in degree $2q-1$ and of finite type
define a central extension of a particular (periodic) thin Lie algebra explicitly constructed in~\cite{CaMa:thin}.
Uniqueness follows upon noting that an infinite-dimensional thin Lie algebra
is centerless, because of the covering property.
Extending this to uniqueness results for other initial patterns in the structure of $(-1)$-algebras
is a crucial step in a prospective classification.
Ideally, one would hope to be able to prove that the earliest genuine diamond of finite type in an infinite-dimensional
$(-1)$-algebra can only occur in certain degrees, and that its degree determines the algebra up to isomorphism.

A uniqueness result for the thin Lie algebra with diamonds of both finite and infinite types
constructed in~\cite{AviMat:-1} and described earlier was proved in~\cite{Avi:thesis},
and is quoted here as Theorem~\ref{thm:presentation}.
However, the hypotheses of Theorem~\ref{thm:presentation}
assume the detailed structure of the thin Lie algebra up to the earliest diamond of finite type,
rather than just its degree.
In particular, this includes that the diamonds of infinite type preceding the earliest genuine diamond of finite type
occur at regular intervals, and that their number is one less than a power of $p$.
The goal of the present paper is a proof that the latter condition is not restrictive.
Put differently, our main result, Theorem~\ref{thm:main}, states that earliest genuine diamond of finite type
of an infinite-dimensional $(-1)$-algebra
can only occur in a degree of the form $(p^{s}+1)(q-1)+1$ for some $s \geq 1$,
under some additional regularity assumptions on the structure of the algebra up to that diamond.
Broadly speaking, the removal of these regularity assumptions is the main obstacle
which separates us from a classification of $(-1)$-algebras, and is the object of current investigation.

Theorem~\ref{thm:main} follows at once from the more precise Theorem~\ref{thm:L-dies},
which includes information on finite-dimensional algebras.
We prove Theorem~\ref{thm:L-dies} in Section~\ref{sec:proof}.
Some of the arguments and calculations introduced in the course of the proof have a wider applicability
and will be used elsewhere.
We briefly comment on this at the end of Section~\ref{sec:proof}.

In our calculations we have been strongly guided by computations performed with
the Australian National University $p$-Quotient program (\cite{HNO}).

\section{General facts on $(-1)$-algebras}\label{sec:preliminaries}%_________________________________________

Throughout the paper $\F$ denotes a field of odd characteristic $p$, with prime field $\F_p$.
We write the Lie bracket in a Lie algebra without a comma
and use the left-normed convention for iterated Lie brackets, so that
$[vyz]=[[vy]z]$.
We also use the shorthand
\[
[v \underbrace{z \ldots z}_{i}]=[vz^i].
\]
The following generalization of the Jacobi identity in a Lie algebra
is easily proved by induction:
\begin{equation}\label{eq:iterated-Jacobi}
[v[y z^{\lambda}]]=\sum_{i=0}^{\lambda}(-1)^i \binom{\lambda}{i}
[vz^iy z^{\lambda-i}].
\end{equation}
This identity will be used repeatedly, usually without specific mention.
Consequently, it will occur very often to evaluate binomial coefficients modulo as elements of $\F$.
This can be done by means of Lucas' Theorem (see~\cite[p.~271]{Dickson1}), which states that
\[
\binom{a}{b} \equiv \prod_{i=0}^{n} \binom{a_i}{b_i} \pmod{p},
\]
where $a=\sum_{i=0}^na_ip^i$ and $b=\sum_{i=0}^nb_ip^i$
are the $p$-adic representations of the non-negative integers $a$ and $b$, thus with $0 \leq a_i,b_i <p$.

We now describe some general features of $(-1)$-algebras and introduce some terminology,
referring the reader to~\cite{CaMa:thin} for further details and the missing proofs.
In particular, we explain how one can assign a {\em type} to each diamond successive to the first.
In this discussion we assume $(-1)$-algebras to be infinite-dimensional unless differently specified,
as this will save us from cumbersome statements.
However, all considerations remain valid for finite-dimensional thin Lie algebras
of dimension large enough; how large should be clear in each instance.
For example, because of the covering property~\eqref{eq:CP} an infinite-dimensional thin Lie algebra is centreless;
however, the centre of a finite-dimensional thin algebra coincides with the
nonzero homogeneous component of highest degree.
We use repeatedly without mention the following more general consequence of the covering property,
for a thin Lie algebra $L$ with $L_1=\langle x,y\rangle$:
whenever an element $u$ spans a one-dimensional component different from the nonzero component of highest degree,
and $[uy]=0$, say, then $[ux]\not=0$.

Let $L= \bigoplus_{i=1}^{\infty} L_i$ be a $(-1)$-algebra with second diamond
in degree $2q-1$, where $q=p^n$ for some $n\geq 1$.
As in the theory of graded Lie algebras of maximal class, a role is played by the centralizers in $L_1$
of homogeneous components, that is,
\[
C_{L_1}(L_k)=\{a \in L_1 | [a,b]=0 \textrm{ for all } b \in L_k\}.
\]
(These are called {\em two-step centralizers} in~\cite{CMN,CN}, a term borrowed from the theory of $p$-groups of maximal class).
Because of the covering property such centralizers are one-dimensional
except for those of diamonds or components immediately preceding a diamond, which are trivial.
It is proved in~\cite{CaJu:quotients} that under our assumptions the one-dimensional components between the first and second
diamond, except that immediately before the latter, all have the same centralizer in $L_1$.
Hence, we may choose $x,y\in L_1$ such that $L_1=\langle x, y \rangle$ and
\begin{equation}\label{eq:two-step}
C_{L_1}(L_2)=C_{L_1}(L_3)= \cdots =C_{L_1}(L_{2q-3})= \langle y \rangle.
\end{equation}

Our assumption that $L$ is infinite-dimensional (or has dimension large enough) implies that the iterated bracket
$v_2=[y x^{2q-3}]$
%\end{equation}
is a nonzero element in $L_{2q-2}$.
The relations
\[
[v_2 xy]+[v_2 yx]=0\quad\textrm{and}\quad
[v_2 yy]=0
\]
hold in $L_{2q}$.
We quote from~\cite{CaMa:thin} the calculations which prove them, as basic
examples of application of the Jacobi identity and its iterate~\eqref{eq:iterated-Jacobi}:
\[
0=[yx^{2q-4}[yxy]]=[yx^{2q-4}[yx]y]-[yx^{2q-4}y[yx]]=
-[yx^{2q-4}xyy]=-[v_2yy],
\]
\begin{align*}
0&=[yx^{q-1}[yx^{q-1}]]=
-\binom{q-1}{q-2}[yx^{2q-3}yx]
+\binom{q-1}{q-1}[yx^{2q-3}xy]\\
&=[v_2 yx]+[v_2 xy].
\end{align*}
These two calculations are simple instances of the two types of calculations which occur many times in this paper.
We use the (iterated) Jacobi relation to expand the Lie bracket $[u,v]$,
where $u$ and $v$ are left-normed iterated Lie brackets in the generators $x$ and $y$ of $L$,
into a left-normed iterated Lie bracket in $x$ and $y$.
In most cases the relation $[u,v]=0$ is obtained as a consequence of a known relation $v=0$,
but sometimes we expand the relation $[u,u]=0$, which is part of the Lie algebra axioms.
We prefer the former type of calculation over the latter whenever there is a choice,
because the former applies in the more general context of {\em Leibniz algebras}
(which are, roughly speaking, Lie algebras without the law $[zz]=0$, see~\cite{Loday:Leibniz}, for example).

Because of the covering property we have $L_{2q}= \langle [v_2 yx] \rangle$,
and so $[v_2 xx]$ is a multiple of $[v_2 yx]$.
If $[v_2 xx]=0$, the second diamond is said to be of {\em infinite type}.
More generally, suppose that $L_h$ is a diamond for some $h>1$,
and suppose that $L_{h-1}$ is one-dimensional, spanned by $w$, say.
The diamond $L_h$ is said to be of \emph{infinite type} if
\begin{equation}\label{eq:type-infinite}
[wxy]+[wyx]=0,\qquad
[wxx]=[wyy]=0.
\end{equation}
The diamond  $L_{h}$ is said to be of \emph{finite type}
 $\lambda $ if
\begin{equation}\label{eq:type-finite}
[wxy]+[wyx]=0,\qquad
[wyy]=0,\qquad
[wyx]= \lambda[wxx],
\end{equation}
for some $\lambda \in \F$.
The latter definition, unlike that of diamond of infinite type, depends on the choice of
$y \in C_{L_1}(L_2)$ and $x \in L_1 \backslash \langle y \rangle$.
A natural normalization, which can be achieved as in~\cite{CaMa:thin} by replacing $y$ with $\lambda^{-1}y$,
is assuming that the first diamond of finite type in order of occurrence (if any) has type $1$.
In the finite-dimensional case, note that when the nonzero homogeneous component of highest degree is two-dimensional
and follows a one-dimensional component,
it formally satisfies the definition of a diamond of any type.

Strictly speaking, a diamond of type zero cannot occur, because then $[wy]$ would be a central element of $L$,
thus violating the covering property.
However, there are situations (here and in other papers, like~\cite{CaMa:thin}, \cite{CaMa:Hamiltonian} and \cite{Avi})
where we have found natural
and convenient to informally call {\em fake diamonds}
certain one-dimensional components $L_h$.
In these cases one has $L_h=\langle w\rangle$ with $[wy]=0$ and $[wxy]=0$.
Therefore, a fake diamond formally satisfies the relations~\eqref{eq:type-finite} with $\lambda=0$,
and will also be referred to as a diamond of type zero.
Note, however, that this last piece of terminology is a matter of convenience,
and not all one-dimensional components satisfying these conditions ought to be dubbed fake diamonds.

We point out a couple of open questions about $(-1)$-algebras for which
we provide partial answers in Section~\ref{sec:chains}.
The first question is whether any diamond in a $(-1)$-algebra can be assigned a type as described above.
In fact, while we call diamond of $L$ any two-dimensional component $L_h$
(besides the fake diamonds described in the previous paragraph),
such a diamond can be assigned a type only if $L_{h-1}$ is one-dimensional, spanned by $w$, say, and the relations
\begin{equation*}
[wxy]+[wyx]=0,\qquad
[wyy]=0
\end{equation*}
hold.
This is likely to always be the case, but we presently have no general proof.
The latter relation is easily proved if there exists $u$ in $L_{h-2}$
such that $[ux]=w$ and $[uy]=0$ (which is typically the case, see the second question below),
because then
\[
0=[u[yxy]]=[u[yx]y]=-[uxyy]=-[wyy].
\]
A proof of the former relation, however, seems to depend on the location and type of the preceding diamond.
We are able to give such a proof only in some specific instances, such as the one
we discuss after Proposition~\ref{prop:chains}.
The second open question is whether one-dimensional components not immediately preceding a diamond
need to be centralized by $y$, in general.
According to~\cite{CaJu:quotients} this is the case for components up to the second diamond, see~\eqref{eq:two-step}.
In Section~\ref{sec:chains} we give a partial positive answer to this question as well.
If both questions will turn out to have positive answers in general,
an arbitrary infinite-dimensional $(-1)$-algebra will be completely described
by specifying the locations and types of all diamonds.

\section{The minimum distance between consecutive diamonds}\label{sec:chains}%_________________________

An argument in the first part of the proof of Proposition~1 in~\cite{CaMa:thin} shows that
a $(-1)$-algebra $T$ cannot have a third (genuine) diamond earlier than in degree $3q-2$.
We simplify that argument and extend it to prove that two consecutive
diamonds in a $(-1)$-algebra must be at least $q-1$ components apart, under certain assumptions
which we have partly discussed at the end of Section~\ref{sec:preliminaries}.

\begin{prop}\label{prop:chains}
Let $T$ be a $(-1)$-algebra in odd characteristic, with second diamond in degree $2q-1$.
Choose generators $x,y\in T_1$ for $T$, with $[yxy]=0$ as usual.
Let  $T_{m}$ be a diamond of $T$ of infinite or nonzero finite type, for some integer  $m\geq 2q-1$.

\begin{enumerate}
\item
If $q>3$ suppose that $y$ centralizes every homogeneous component from $T_{m-q+2}$ up to $T_{m-2}$.
Then $y$ centralizes every homogeneous component from $T_{m+1}$ up to $T_{m+q-3}$.
In particular, all components $T_{m+1},\ldots,T_{m+q-2}$ are at most one-dimensional,
and so a diamond after $T_m$ cannot occur in degree lower than $m+q-1$.
\item
If $T_{m+q-2}=\langle w\rangle$ then the relations
$[wxy]+[wyx]=0$ and $[wyy]=0$ hold.
If $[wy]\neq 0$ and $T_{m+q}\neq\{0\}$ then $T_{m+q-1}$ is a genuine diamond,
and can be assigned a type.
\end{enumerate}
\end{prop}

\begin{proof}
Let $v$ be a nonzero element in $T_{m-1}$.
Since $T_{m}$ is a diamond of infinite or finite nonzero type, say $\lambda$, the relations
\[
[vxy]+[vyx]=0,\qquad [vyy]=0,\qquad
\lambda^{-1}[vyx]=[vxx]
\]
hold in degree $m+1$, where we read $\infty^{-1}=0$.
The actual value of $\lambda$ will play no role in most of the proof.
Recall from Section~\ref{sec:preliminaries} that in a $(-1)$-algebra we have
$[yx^{i}y]=0$ for $1\le i\le 2q-4$.
Differently from the proof of Proposition~1 in \cite{CaMa:thin}, which used several of these relations,
here we use the single relation
$[yx^{q-2}y]=0$
to prove by induction on $h$ that
\[
[vyx^hy]=0 \quad \textrm{for} \quad 1 \leq h \leq q-3.
\]
This will imply, inductively, that $T_{m+h}$ is spanned by $[vyx^h]$,
for $1\le h\le q-3$, and that it is centralized by $y$,
thus proving assertion (1).

We prove the induction base and step at the same time.
Let $1\le h\le q-3$, and in case $h>1$ suppose that the conclusion holds up to $h-1$.
We may choose $u$ such that $v=[ux^{q-2-h}]$.
Note that $q-2-h>0$ and that $[ux^iy]=0$ for $0\le i<q-2-h$.
Expanding
\begin{align*}
0&=[u[yx^{q-2}y]]= [u[yx^{q-2}]y]\\
&=(-1)^{h+1}\binom{q-2}{h}[vyx^hy]+(-1)^{h}\binom{q-2}{h-1}[vxyx^{h-1}y]\\
&=(-1)^{h+1}\binom{q-1}{h}[vyx^hy]=-[vyx^hy]
\end{align*}
we reach the desired conclusion that $[vyx^hy]=0$.

To prove assertion (2), note that $T_{m+q-2}=\langle w\rangle$, where $w=[vxyx^{q-3}]$.
The computations
\[
0=[vyx^{q-3}[yxy]]=-[vyx^{q-2}yy]=[vxyx^{q-3}yy]=[wyy]
\]
and
\begin{align*}
0&=[vx[yx^{q-2}y]]= [vx[yx^{q-2}]y]-[vxy[yx^{q-2}]]\\
&=\binom{q-2}{0}[vxyx^{q-2}y]-\binom{q-2}{q-2}[vxxx^{q-3}yy]\\
&\quad
-\binom{q-2}{q-3}[vxyx^{q-3}yx]+\binom{q-2}{q-2}[vxyx^{q-2}y]\\
&=2[wxy]+2[wyx]+\lambda^{-1}[wyy]
\end{align*}
produce the desired relations.
Now suppose that $[wy]$ spans $T_{m+q-1}$,
whence $[wx]=\nu[wy]$ for some scalar $\nu$.
Then $[wyy]=0$ and
\[
0=[wxy]+[wyx]=\nu[wyy]+[wyx]=[wyx]
\]
imply that $[wy]$ is a central element of $T$, and so $T_{m+q}=\{0\}$ by the covering property.
Therefore, if
$[wy]\neq 0$ and $T_{m+q}\neq\{0\}$
then $T_{m+q-1}$ is two-dimensional, and because of the relations proved above it can be assigned a type.
\end{proof}

\begin{rem}
We have excluded the case $q=3$ in assertion~(1) of Proposition~\ref{prop:chains}
because then both its hypothesis and conclusion would be void, apart from the final consequence
that a diamond after $T_m$ cannot occur in degree lower than $m+q-1$.
This last assertion holds when $q=3$ as well, because assuming that $T_m$ has a type implies that $T_{m+1}$
is one-dimensional.

Note that when $q=3$ the $(-1)$-algebra $T$ of Proposition~\ref{prop:chains}
has second diamond in degree $5$.
However, the uniqueness result from~\cite{CMNS} for such algebras, which we have
mentioned in the Introduction, only holds in characteristic greater than five.
In fact, several infinite families of pairwise non-isomorphic $(-1)$-algebras
with $q=3$ arise as special cases of various constructions~\cite{CaMa:Hamiltonian,AviMat:-1}.
\end{rem}

Proposition~\ref{prop:chains} can be used to prove, inductively, that all one-dimensional homogeneous components of $T$
not immediately preceding a diamond are centralized by $y$, and that each diamond can be assigned a type,
thus giving positive answers to both questions discussed at the end of Section~\ref{sec:preliminaries},
but only as far as all the diamonds involved (from the second on) occur at regular distance of $q-1$ components apart.
However, consecutive diamonds of a $(-1)$-algebra can be farther apart
(in which case $[wy]=0$ holds in the last part of Proposition~\ref{prop:chains}),
and this breaks the inductive argument outlined above.
It is known, from the connection with graded Lie algebras of maximal class established in~\cite{CaMa:thin},
that the difference in the degrees of consecutive diamonds of a $(-1)$-algebra can take all positive values of the form
$2q-1-p^r$, with $r$ a nonnegative integer or $-\infty$.
It seems likely that these are all possibilities for the distance between consecutive diamonds,
but we are presently unable to produce a proof.

\section{$(-1)$-algebras with diamonds of finite and infinite type}\label{sec:mixed}%_________________________________________

In this section we describe a particular family of infinite-dimensional $(-1)$-algebras
having diamonds of both finite and infinite type
which was identified in~\cite{Avi:thesis}.
As we said in the Introduction, we believe that all infinite-dimensional $(-1)$-algebras
having diamonds of both finite and infinite type should belong to this family,
and our Theorem~\ref{thm:main} is a step in this direction.

Let $L$ be an infinite-dimensional $(-1)$-algebra with diamonds of infinite type in all degrees
$k(q-1)+1$ for $k=2, \ldots, p^s$, where $s \geq 1$,
and with a diamond of finite type $\lambda\in \F^*$ in degree $(p^s+1)(q-1)+1$.
By replacing $y$ with $\lambda^{-1}y$ we may assume that $\lambda=1$.
According to~\cite{Avi:thesis}, the algebra $L$ turns out to be uniquely determined by these prescriptions.
It has diamonds in all degrees of the form $t(q-1)+1$.
If $t \not \equiv 1 \pmod{p^s}$ the corresponding diamond is of infinite type.
If $t \equiv 1 \pmod{p^s}$, say $t= rp^s+1$, the corresponding diamond
has type $r$ (viewed modulo $p$), hence an element of the prime field $\F_p$.
In particular, when $r \equiv 0 \pmod{p}$ we have a fake diamond.
The diamonds of finite type (including the fake ones, which have type zero) follow an arithmetic progression.
Furthermore, each one-dimensional homogeneous component not immediately preceding a genuine diamond
is centralized by $y$.
This uniqueness result was proved in \cite{Avi:thesis} by showing that if an algebra $N$
is defined by a finite presentation encoding part of the above prescriptions up to degree $(p^s+1)(q-1)+2$,
then the quotient $L$ of $N$ modulo its centre is a thin algebra and has the structure stated above.
We quote the precise result from~\cite{Avi:thesis}.

\begin{theorem}\label{thm:presentation}
Let $q$ be a power of the odd prime $p$.
Let $N=\displaystyle{\bigoplus_{i=1}^{\infty}N_{i}}$ be the Lie algebra on two generators $x$ and $y$
subject to the following relations, and graded by assignign degree $1$ to $x$ and $y$:
\begin{align*}
&[y x^{i}y]=0&  &\textrm{for $0 < i <2q-3$ with $i\not= 2q-p^t-2$,} & & \\
&[y x^{2q-p^t-2}yx]=0&      &\textrm{for $1 \leq t \leq n$,} & & \\
&[yxyy]=0             & & \textrm{if $q=p=3$,}  \\
&[v_2xxx]=0=[v_2xxy],&    &  & &\\
&[v_kxx]=0&    &\textrm{for $3 \leq k \leq p^s$ with $k$ even,} & &\\
&[v_{p^s+1}yx]=[v_{p^s+1}xx].& & &
\end{align*}
Here $v_k$ is defined recursively by
$v_2=[y x^{2q-3}]$, and
$v_{k}=[v_{k-1} xy x^{q-3}]$ for $k>2$.
Then $L=N/\cent(N)$ is a thin algebra and has the diamond structure described above in the text.
In particular, a thin algebra generated by the homogeneous elements $x$ and $y$
and satisfying the above relations is necessarily isomorphic with $L$.
\end{theorem}

Note that the elements $v_k$ have degree $k(q-1)$
and span the components just preceding the (possibly fake) diamonds.
The last three relations given in Theorem~\ref{thm:presentation}
specify the types of the diamonds in the degrees $k(q-1)+1$, for $k=2, \ldots, p^s+1$.

The Lie algebra $N$ of Theorem~\ref{thm:presentation} is constructed in~\cite{AviMat:-1}
as a loop algebra of a certain finite-dimensional Lie algebra.
It is also proved there that its central quotient $L$ cannot be finitely presented.
In particular, it is not possible to supplement the presentation given in Theorem~\ref{thm:presentation}
by adding finitely many relations and obtain a presentation of the thin Lie algebra $L$.
In the opposite direction, some of the relations in the presentation turn out to be superfluous,
as we discuss in Section~\ref{sec:first-finite}.

In view of the preceding discussion the last claim of Theorem~\ref{thm:presentation} can be rephrased as follows:
up to isomorphism $L=N/\cent(N)$ is the unique infinite-dimensional thin algebra
with second diamond in degree $2q-1$, diamonds of infinite type in degrees
$i(q-1)+1$, for $2 \leq i \leq p^s$, and a diamond of type one
in degree $(p^s+1)(q-1)+1$.

\section{The degree of the first diamond of finite type}\label{sec:first-finite}%_________________________

In this section we state our main results.
Consider a $(-1)$-algebra $T$ with diamonds of infinite type in each degree $k(q-1)+1$ for $k=2, \ldots, a-1$,
and  a diamond in degree $a(q-1)+1$, for some integer $a \geq 3$.
We determine for which values of $a$ the diamond in degree
$a(q-1)+1$ can be prescribed to be of nonzero finite type
(and, hence, of type 1, possibly after adjusting the generators of $T$ as explained in~Section \ref{sec:preliminaries})
without the algebra being forced to have finite dimension.
More precisely, we prove the following result.

\begin{theorem}\label{thm:L-dies}
Let $T$ be a $(-1)$-algebra in odd characteristic with second diamond in degree $2q-1$
and with $T_{a(q-1)+2}\neq 0$, where $a\geq 3$.
Suppose that $T$ has diamonds of infinite type in each degree $k(q-1)+1$ for $k=2, \ldots, a-1$,
and a diamond of nonzero finite type in degree $a(q-1)+1$.
Then the following assertions hold:
\begin{enumerate}
\item $a$ is even;
\item if $a \not \equiv 1 \bmod p$ then $T_{(a+1)(q-1)+3}=0$;
\item if $a \equiv 1 \bmod p$, but $a-1$ is not a power of $p$, then $T_{(a+p^s)(q-1)+2}=0$,
where $p^s$ is the highest power of $p$ which divides $a-1$.
\end{enumerate}
\end{theorem}

Assertions (2) and (3) of Theorem~\ref{thm:L-dies} show that $T$ has bounded dimension unless $a-1$ is a power of $p$.
Assertion (1) will play a crucial role in the proof of assertion (2) given in Lemma~\ref{lemma:a-notcong-1}.
As we have mentioned near the end of Section~\ref{sec:first-finite},
it also justifies the omission of the relations $[v_kxx]=0$ for $k$ odd
in the presentation for a central extension of $T$ given in~Theorem~\ref{thm:presentation}.

\begin{rem}
We point out an additional piece of information which is not explicit in Theorem~\ref{thm:L-dies}.
When $a$ is odd, the assumptions of Theorem~\ref{thm:L-dies} can all be satisfied except for the requirement that
the diamond in degree $a(q-1)+1$ has nonzero finite type.
Assuming only that $T$ has a genuine diamond in that degree, it follows that such diamond must have infinite type.
This is essentially the content of Lemma~\ref{lemma:a-odd}, but also follows from the present formulation:
that diamond can be assigned a type because of Proposition~\ref{prop:chains},
and the type must be infinite because of Theorem~\ref{thm:L-dies}.
\end{rem}

The special case of Theorem~\ref{thm:main} where $T$ is infinite-dimensional
reduces to the following more concise statement.

\begin{theorem}\label{thm:main}
The first diamond of nonzero finite type (if any) of an infinite-dimensional $(-1)$-algebra $T$
in odd characteristic occurs in degree  $(p^{s}+1)(q-1)+1$, for some $s\geq 0$, provided
$T$ has diamonds of infinite type in all lower degrees congruent to $1$ modulo $q-1$, except $q$.
\end{theorem}

As we have discussed in Section~\ref{sec:mixed},
under the assumptions of Theorem~\ref{thm:main} the $(-1)$-algebra $T$ is then isomorphic to one of the algebras $L$
of Theorem~\ref{thm:presentation}.

There is some redundancy in the hypotheses of Theorem~\ref{thm:L-dies}
on the structure of $T$ up to the diamond of finite type,
and we point out here what needs to be assumed and what is a consequence.
This is of interest when writing an efficient presentation for $T$, or possibly some algebra which has $T$ as a graded quotient
(as in Theorem~\ref{thm:presentation}),
which was an essential task in carrying out the computational experiments which lead to the present results.
In particular, this will also shed some light on the presentation given in Theorem~\ref{thm:presentation}
(for which a complete proof is given in~\cite{Avi:thesis}).
At the same time we fix some notation which we will use in the proof of Theorem~\ref{thm:L-dies} in the next section.

Let $T$ be a $(-1)$-algebra, with second diamond in degree $2q-1$,
and choose  $x,y \in T_{1}$ such that $T_{1}=\langle x, y \rangle$ and $[yxy]=0$.
As we have already mentioned in Section~\ref{sec:preliminaries} in the infinite-dimensional case,
it follows from~\cite{CaJu:quotients} that $\langle y \rangle =C_{T_{1}}(T_{i})$ for  $2\leq i \leq 2q-3$,
provided $T_{3q-1}\not=\{0\}$.
We assume the dimension (or, equivalently, the nilpotency class)
of $T$ to be large enough throughout the argument and set
\[
v_{1}=[yx^{q-2}]\quad\textrm{and}\quad v_{2}=[v_{1}x^{q-1}]=[yx^{2q-3}].
\]
Thus, $v_{1}$ and $v_{2}$ are nonzero elements in degree $q-1$ and $2q-2$.
We have seen in Section~\ref{sec:preliminaries} that
\[
[v_{2}xy]+[v_{2}yx]=0\quad\textrm{and}\quad [v_{2}yy]=0.
\]
We specify the diamond $T_{2q-1}$ to have infinite type
by imposing that $[v_{2}xx]=0$.

Define recursively the elements
\begin{equation}\label{eq:vk+1}
v_{k}=[v_{k-1} xy x^{q-3}] \quad\textrm{for}\quad  3 \leq k \leq a.
\end{equation}
Thus, $v_k$ has degree $k(q-1)$.
We now show how a few carefully chosen relations involving the elements $v_k$
suffice to determine the structure of $T$ up to degree $a(q-1)+2$,
satisfying the hypotheses of Theorem~\ref{thm:L-dies}.
Let $3\le k\le a$.
Assume recursively that
\[
[v_{k-2}xyx^{h}y]=0 \quad\textrm{for}\quad 0\leq h \leq q-4
\]
(but $[v_{1}xxx^{h}y]=0$ instead in case $k=3$),
whence $v_{k-1}$ spans $T_{(k-1)(q-1)}$ by the covering property.
Assume also that $T_{(k-1)(q-1)+1}$ is a diamond of infinite type.
Then assertion (1) of Proposition~\ref{prop:chains} yields that
\[
[v_{k-1}xyx^{h}y]=0 \quad\textrm{for}\quad 0\leq h \leq q-4.
\]
In particular, $v_{k}$ spans $T_{k(q-1)}$ because of the covering property, and hence
$T_{k(q-1)+1}$ is spanned by $[v_{k}x]$ and $[v_{k}y]$.
Assertion (2) of Proposition~\ref{prop:chains} then implies that
\[
[v_{k}xy]+[v_{k}yx]=0\quad\textrm{and}\quad [v_{k}yy]=0.
\]
In order to complete one step of our recursive definition of $T$ we need only
impose the single relation $[v_{k}xx]=0$ if $k<a$, and $\lambda^{-1}[v_{a}yx]=[v_{a}xx]$ otherwise.
Our assumption that the dimension of $T$ is large enough then implies that $[v_{k}y]\neq 0$,
otherwise the relations would yield
$T_{k(q-1)+2}=\langle[v_{k}xx],[v_{k}xy]\rangle=\{0\}$.
Then $T_{k(q-1)+1}$ is a diamond, again according to assertion (2) of Proposition~\ref{prop:chains},
and it has the type which we have specified.

We conclude that if a thin algebra $T$ with second diamond in degree $2q-1$
satisfies $[v_kxx]=0$ for $2\le k<a$ and $\lambda^{-1}[v_ayx]=[v_axx]$, where the elements $v_k$ are defined as above, then
$T$ has diamonds in all degrees $k(q-1)+1$ for $2\le k\le a$, all of infinite type except the last, perhaps,
which has type $\lambda$.
Furthermore, all one-dimensional components of $T$ which do not immediately precede a diamond are centralized by $y$.
It will follow from Lemma~\ref{lemma:a-odd} that the relations $[v_kxx]=0$ for $k$ odd are actually superfluous
and, in fact, they are omitted in the presentation for a central extension of $T$ given in~Theorem~\ref{thm:presentation}.

It is worth anticipating here that if the recursive definition of the elements
$v_k$ in~\eqref{eq:vk+1} is extended past $v_{a}$, the inductive argument given above carries through
for the subsequent diamonds and chains on the only assumption that for each $k$ some relation of the form
$\mu^{-1}[v_kyx]=[v_kxx]$, with $\mu$ finite nonzero or infinity, is imposed or proved by other means.
This is because the specific type of each diamond $T_{k(q-1)+1}$ was immaterial in the argument.
We will use this extension in the proof of Lemma~\ref{lemma:a-cong-1}.

\section{Proof of Theorem~\ref{thm:L-dies}}\label{sec:proof}%_________________________

We divide the proof into several lemmas.
After the following technical lemma describing the adjoint action of $v_1$ and $v_2$
on homogeneous elements close to diamonds,
we prove assertions (1), (2) and (3) of Theorem~\ref{thm:L-dies} in
Lemmas~\ref{lemma:a-odd}, \ref{lemma:a-notcong-1} and~\ref{lemma:a-cong-1}, respectively.

\begin{lemma}\label{lemma:technical}
Let the notation and assumptions be as in Section~\ref{sec:first-finite}, and let $k\ge 2$.
\begin{enumerate}
\item
Suppose that $[v_kyx]=\mu[v_kxx]$, with $\mu$ finite nonzero, or infinity, that is,
the diamond $T_{k(q-1)+1}$ has type $\mu$.
Then we have
\begin{align*}
&[v_kv_1]=v_{k+1},\\
&[v_kxv_1]=[v_{k+1}x]+\mu^{-1}[v_{k+1}y],&
&[v_kyv_1]=[v_{k+1}y],\\
&[v_kxxv_1]=\mu^{-1}[v_{k+1}yx],&
&[v_kyxv_1]=[v_{k+1}yx].&
\end{align*}
\item
Suppose that $[v_kxx]=0$ and $[v_{k+1}xx]=0$, that is, both diamonds $T_{k(q-1)+1}$ and $T_{(k+1)(q-1)+1}$
have infinite type.
Then we have
\[
[v_kv_2]=0,\quad [v_kxv_2]=0,\quad [v_kyv_2]=0,\quad
[v_kyxv_2]=[v_{k+2}xx],
\]
and if $q>3$ also
\[
[v_kxyx^{q-4}v_1]=[v_{k+1}xyx^{q-4}],\quad
[v_kxyx^{q-4}v_2]=3[v_{k+2}x^{q-2}].
\]
\end{enumerate}
\end{lemma}

\begin{proof}
The first of the set of equalities in (1) is
\[
[v_kv_1]=[v_k[yx^{q-2}]]
=[v_k yx^{q-2}]+2[v_k xyx^{q-3}]
=v_{k+1},
\]
and the remaining ones can be proved similarly.

To prove (2) suppose that $[v_kxx]=0$ and $[v_{k+1}xx]=0$.
We have
\[
[v_kv_2]=[v_k[yx^{2q-3}]]=[v_kyx^{2q-3}]+3[v_kxyx^{2q-4}]=2[v_{k+1}x^{q-1}]=0,
\]
and
\[
[v_kyv_2]=
\binom{2q-3}{q-2}[v_{k+1}yxx^{q-2}]-\binom{2q-3}{q-1}[v_{k+1}xyx^{q-2}]=0,
\]
because of the vanishing modulo $p$ of the binomial coefficients involved.
The other two equalities are similar.
Finally, if $q>3$ we have
$$
[v_kxyx^{q-4}v_1]
=-\binom{q-2}{1}[v_{k+1}yxx^{q-4}]+\binom{q-2}{2}[v_{k+1}xyx^{q-4}]
=[v_{k+1}xyx^{q-4}],
$$
and
$$
[v_kxyx^{q-4}v_2]
=-\binom{2q-3}{1}[v_{k+1}yxx^{2q-5}]+\binom{2q-3}{2}[v_{k+1}xyx^{2q-5}]
=3[v_{k+2}x^{q-2}],
$$
which completes the proof.
\end{proof}
According to the very first equality proved in Lemma~\ref{lemma:technical},
the recursive definition~\eqref{eq:vk+1} of $v_k$ can now be replaced by the more compact formula
$v_k=[v_2v_1^{k-2}]$.
We can use this, together with the other equalities given in Lemma~\ref{lemma:technical}, to
compute the adjoint action of any $v_h$ on homogeneous elements close to diamonds.
One instance of this type of calculation occurs in the last part of the proof of the following lemma.
More instances will occur in the rest of this section,
and will be treated in less detail.

\begin{lemma}\label{lemma:a-odd}
Under the hypotheses of Theorem~\ref{thm:L-dies}, if $a$ is odd then the relation
$[v_axx]=0$ holds in degree $a(q-1)+2$.
\end{lemma}
\begin{proof}
We expand the relation $0=[u,u]$, where $u$ is a homogeneous element
of degree $a(q-1)/2+1$.
We distinguish two cases.
For the case $a=3$ we choose $u=[v_{1}x^{\frac{q+1}{2}}]=
[yx^{\frac{3q-3}{2}}]$
and we have
$$
0=[u[yx^{\frac{3q-3}{2}}]]=\sum_{i=0}^{\frac{3q-3}{2}}(-1)^{i}
\binom{\frac{3q-3}{2}}{i}[ux^{i}yx^{\frac{3q-3}{2}-i}]
=\pm \binom{\frac{3q-3}{2}}{\frac{q-3}{2}}
[v_2yx^q]=\pm [v_3xx].
$$
Now suppose that $a>3$ and write $a=2h+1$. We choose
$u=[v_{h}xyx^{\frac{q-3}{2}}]$ and we have
\begin{align*}
0&=[u[v_{h}xyx^{\frac{q-3}{2}}]]=\sum_{i=0}^{\frac{q-3}{2}}(-1)^i
\binom{\frac{q-3}{2}}{i}
[ux^{i}[v_hxy]x^{\frac{q-3}{2}-i}]& &  \\
&= \sum_{i=0}^{\frac{q-3}{2}}(-1)^i
\binom{\frac{q-3}{2}}{i}
\left([ux^{i}[v_hx]yx^{\frac{q-3}{2}-i}]
-[ux^{i}y[v_hx]x^{\frac{q-3}{2}-i}]\right).
\end{align*}
The element $[ux^{i}[v_hx]]$ belongs to $T_{2h(q-1)+\frac{q+3}{2}+i}$,
and so is centralized by $y$ unless $i=\frac{q-5}{2}, \frac{q-3}{2}$.
Similarly, the element $[ux^{i}y]$ vanishes except, perhaps,
when $i=\frac{q-3}{2}$.
Therefore, we have
\begin{align*}
0&=(-1)^{\frac{q-5}{2}}[u,u]
= \binom{\frac{q-3}{2}}{\frac{q-5}{2}}
[v_hxyx^{q-4}[v_hx]yx]
- [v_{h+1}[v_hx]y]+ [v_{h+1}y[v_hx]]\\
&= \binom{\frac{q-3}{2}}{\frac{q-5}{2}}[v_hxyx^{q-4}v_hxyx]
- \binom{\frac{q-3}{2}}{\frac{q-5}{2}}[v_{h+1}v_hyx]+\\&\quad
-[v_{h+1}v_hxy]+[v_{h+1}xv_hy]
+[v_{h+1}yv_hx]-[v_{h+1}yxv_h]\\
&=-[v_{h+1}yxv_h]=(-1)^{h+1}[v_axx],
\end{align*}
which gives the desired conclusion.
The final steps of the above calculation depend on the following relations, which we prove using Lemma~\ref{lemma:technical}
several times.
If $q>3$ we have
\begin{align*}
[v_hxyx^{q-4}v_hxyx]&=[v_hxyx^{q-4}[v_2{v_1}^{h-2}]xyx]\\
&=\sum_{i=0}^{h-2}(-1)^i\binom{h-2}{i}[v_{h+i}xyx^{q-4}v_2{v_1}^{h-2-i}xyx]\\
&=3\sum_{i=0}^{h-2}(-1)^i\binom{h-2}{i}[v_{h+i+2}x^{q-2}{v_1}^{h-2-i}xyx]=0,
\end{align*}
since $h+i+2 \leq 2h$ and $[v_{h+i+2}xx]=0$.
If $q=3$ this term may be interpreted as $-[v_hyv_hxyx]$,
but it does not really appear in the earlier calculation, because
its coefficient $\binom{(q-3)/2}{(q-5)/2}$ vanishes in this case.
Similar computations show that
\[
[v_{h+1}v_hyx]=0=[v_{h+1}xv_hy]=[v_{h+1}yv_hx].
\]
Finally, the calculation
\begin{align*}
[v_{h+1}yxv_h]&=\sum_{i=0}^{h-2}(-1)^i\binom{h-2}{i}[v_{h+i+1}yxv_2{v_1}^{h-2-i}]\\
&=\sum_{i=0}^{h-2}(-1)^i\binom{h-2}{i}[v_{h+i+3}x^2{v_1}^{h-2-i}]\\
&=(-1)^{h-2}[v_{a}xx]
\end{align*}
completes the proof.
\end{proof}

The relation $[v_axx]=0$ proved in Lemma~\ref{lemma:a-odd} contradicts our assumption that the diamond
in degree $a(q-1)+1$ has finite nonzero type.
Hence $a$ must be even, and assertion~(1) of Theorem~\ref{thm:L-dies} is proved.

Now assume that $T$ has a diamond of nonzero finite type in degree $a(q-1)+1$, with $a>2$,
and that $T$ has diamonds of infinite type
in all lower degrees congruent to $1$ modulo $q-1$, with the exception of $q$.
Assume also that $T_{(a+1)(q-1)+3}\neq \{0\}$.
According to Lemma~\ref{lemma:a-odd}, $a$ must be even.
As we may, we assume that the type of the finite diamond is $1$, and so the relation
\[
[v_axx]=[v_ayx]
\]
holds.
The discussion after the statement of Theorem~\ref{thm:L-dies} shows that
$[v_ayy]=0=[v_axy]+[v_ayx]$, that
\[
[v_ayx^hy]=0 \quad\text{for}\quad 1 \leq h \leq q-3
\]
and, finally, that
$[v_{a+1}yy]=0=[v_{a+1}xy]+[v_{a+1}yx]$.
These relations imply that the (possibly fake) diamond $T_{(a+1)(q-1)+1}$ can be assigned a type.
Now we show that this type can only be infinity (and, in particular, the diamond is not fake).
This follows from the relations
\[
0=[v_{a-1}x[v_2xx]]=[v_{a-1}xv_2xx]
=[v_{a-1}xyx^{2q-1}]=-[v_{a+1}xxx]
\]
and
\begin{align*}
0&=[v_{a-1}y[v_{2}xx]]=[v_{a-1}y[yx^{2q-1}]]\\
&=[v_{a-1}yx^{q-2}yxx^{q}]+[v_{a-1}yx^{q-2}xyx^{q}]+\\
&\quad-[v_{a-1}yx^{2q-3}yxx]-[v_{a-1}yx^{2q-3}xyx]-[v_{a-1}yx^{2q-3}xxy]\\
&=-[v_ayxx^q]-[v_axyx^q]+[v_ax^{q-1}yxx]+[v_ax^{q-1}xyx]+[v_ax^{q-1}xxy]\\
&=-[v_{a+1}yxx]-[v_{a+1}xyx]-[v_{a+1}xxy]\\
&=-[v_{a+1}xxy].
\end{align*}
In fact, these relations imply that $[v_{a+1}xx]$ is central, and so
$[v_{a+1}xx]=0$ because of the covering property.
It also follows that $[v_{a+1}x]$ and $[v_{a+1}y]$ are not proportional, otherwise they would be central,
and hence vanish, contradicting our assumption that $T_{(a+1)(q-1)+3}\neq\{0\}$.
We conclude that $T_{(a+1)(q-1)+1}$ is a diamond of infinite type, as claimed.

The following lemma proves assertion~(2) of Theorem~\ref{thm:L-dies}.
\begin{lemma}\label{lemma:a-notcong-1}
Under the hypotheses of Theorem~\ref{thm:L-dies}, if $a \not \equiv 1 \pmod{p}$
and the diamond in degree
$a(q-1)+1$ is of nonzero finite type,
then $T_{(a+1)(q-1)+3}=\{0\}$.
\end{lemma}

\begin{proof}
We work by way of contradiction and assume that $T_{(a+1)(q-1)+3}\neq\{0\}$.
We compute $[v_{a}xxv_1]$ in two different ways starting from the obvious relation $[v_{a}xxv_1]=-[v_1[v_{a}xx]]$.
According to Lemma~\ref{lemma:technical}
we have
\[
[v_{a+1}yx]=[v_{a}xxv_1]=-[v_1[v_{a}xx]]=2[v_1xv_ax]-[v_1xxv_a],
\]
since $[v_1v_{a}xx]=-[v_{a+1}xx]=0$.
The two Lie brackets at the right-hand side become
\begin{align*}
[v_1xv_ax]&=[v_1x[v_2v_1^{a-2}]x]
=-\sum_{i=1}^{a-2}(-1)^i\binom{a-2}{i}[v_2yv_1^{i-1}v_2v_1^{a-2-i}x]\\
&=-[v_{a-1}yv_2x]=[v_{a+1}yx]
\end{align*}
and
\begin{align*}
[v_1xxv_a]&=[v_1xx[v_2v_1^{a-2}]]
=-\sum_{i=1}^{a-2}(-1)^i\binom{a-2}{i}[v_2yxv_1^{i-1}v_2v_1^{a-2-i}]\\
&=(a-2)[v_2yxv_1^{a-4}v_2v_1]-[v_2yxv_1^{a-3}v_2]\\
&=(a-2)[v_{a-2}yxv_2v_1]-[v_{a-1}yxv_2]\\
&=(a-2)[v_{a}xxv_1]+2[v_{a+1}yx]
=a[v_{a+1}yx].
\end{align*}
Besides the fact that $a$ is even, here we have used repeatedly Lemma~\ref{lemma:technical}, supplemented by the relations
\begin{align*}
&[v_1xv_1]=-[v_2y],&
&[v_1xxv_1]=-[v_2yx],\\
&[v_1xv_2]=0,&
&[v_1xxv_2]=0,\\
&[v_{a-1}yv_2]=-[v_{a+1}y],&
&[v_{a-1}yxv_2]=-2[v_{a+1}yx],
\end{align*}
which can be verified similarly.
We conclude that
$[v_{a+1}yx]=(2-a)[v_{a+1}yx]$,
whence $a\equiv 1\pmod{p}$.
This contradicts one of our hypotheses.
\end{proof}

The following lemma completes the proof of Theorem~\ref{thm:L-dies} by establishing assertion~(3).

\begin{lemma}\label{lemma:a-cong-1}
Under the hypotheses of Theorem~\ref{thm:L-dies}, suppose that $a \equiv 1 \pmod{p}$ and
that the diamond in degree
$a(q-1)+1$ is of type $1$.
Write $a=1+np^s$, for some integer $s \geq 1$, with
$n \not \equiv 0 \pmod{p}$.
If $n>1$ then $T_{(a+p^s)(q-1)+2}=\{0\}$.
\end{lemma}

\begin{proof}
This proof will be in a terser style than the previous ones, as we have encountered several similar
calculations before.
In particular, see the comment which precedes Lemma~\ref{lemma:a-odd}.

Extending the recursive definition~\eqref{eq:vk+1} we set
\[
v_{a+k}=[v_{a+k-1}xyx^{q-3}] \quad\text{for}\quad 2 \leq k \leq p^s.
\]
As we have mentioned at the end of the discussion which follows the statement of Theorem~\ref{thm:L-dies},
the inductive argument described there
can be extended past $v_{a+1}$ to prove that for $k=1,\ldots,p^s$ we have
\begin{align*}
&[v_{a+k}yy]=0=[v_{a+k}xy]+[v_{a+k}yx],\quad\text{and}\\
&[v_{a+k}yx^h y]=0 \quad\text{for}\quad 1\leq h\leq q-3,
\end{align*}
provided we can show that each component $\langle[v_{a+k}x],[v_{a+k}y]\rangle$ is actually two-dimensional,
that is, a diamond.
This follows from the relations
\begin{equation}\label{eq:va+kxx}
[v_{a+k}xx]=0 \quad\text{for}\quad 1 \leq k \leq p^s,
\end{equation}
which we prove now (and also show that those diamonds have infinite type).
As usual, we do this inductively, assuming that all claimed relations in lower degrees hold.
The case $k=1$ was proved before Lemma~\ref{lemma:a-notcong-1}.

Since $a=np^s+1>p^s+1$ we have $[v_{k+1}xx]=0$ for $k=2, \ldots p^s$.
Consequently, we have
\begin{align*}
0&=[v_{a-1}[v_{k+1}xx]]
=[v_{a-1}[v_2v_1^{k-1}]xx]-2[v_{a-1}x[v_2v_1^{k-1}]x]\\
&=[v_{a-1}v_2v_1^{k-1}xx]-(k-1)[v_{a}v_2v_1^{k-2}xx]\\
&\quad
-2[v_{a-1}xv_2v_1^{k-1}x]+2(k-1)[v_{a}xv_2v_1^{k-2}x]\\
&=(1-k)[v_{a+k}xx].
\end{align*}
Besides the inductive hypotheses we have used Lemma~\ref{lemma:technical} repeatedly,
supplemented in the last step by the additional relations
\begin{align}\label{eq:va-1}
&[v_{a-1}v_2]=-2\,v_{a+1},&
&[v_{a-1}xv_2]=-[v_{a+1}x],\\ \label{eq:va}
&[v_{a}v_2]=v_{a+2},&
&[v_{a}xv_2]=0,
\end{align}
which can be verified similarly.
We conclude that $[v_{a+k}xx]=0$ as desired, except, perhaps, when $k \equiv 1 \pmod{p}$.
To deal with these exceptions we write $k=hp^t+1$
with $h \not \equiv 0 \pmod{p}$.
Since $k+p^t<2k\le 2p^s< a$
we have
\begin{align*}
0&=[v_{a-p^t}[v_{k+p^t}xx]]
=[v_{a-p^t}[v_2v_1^{k+p^t-2}]xx]-2[v_{a-p^t}x[v_2v_1^{k+p^t-2}]x]\\
&=\binom{k+p^t-2}{p^t-1}[v_{a-1}v_2v_1^{k-1}xx]-\binom{k+p^t-2}{p^t}[v_{a}v_2v_1^{k-2}xx]\\
&\quad
-2\binom{k+p^t-2}{p^t-1}[v_{a-1}xv_2v_1^{k-1}x]+2\binom{k+p^t-2}{p^t}[v_{a}xv_2v_1^{k-2}x]\\
&=-h[v_{a+k}xx].
\end{align*}
and this completes a proof of~\eqref{eq:va+kxx}.

In particular, we have proved that
\[
[v_{a+p^s}xy]+[v_{a+p^s}yx]=
[v_{a+p^s}yy]=[v_{a+p^s}xx]=0,
\]
and hence the component of $T$ of degree $(a+p^s)(q-1)+2$ is spanned by $[v_{a+p^s}xy]$.
To complete the proof we only need to show that this element vanishes, and we do so by expanding
both sides of the equality
\[
[v_{p^s}[v_ayx]]=
[v_{p^s}[v_axx]].
\]
Using Lemma~\ref{lemma:technical} and relations~\eqref{eq:va-1} and~\eqref{eq:va} one sees
that the right-hand side,
\[
[v_{p^s}[v_axx]]=[v_{p^s}[v_2v_1^{a-2}]xx]-2[v_{p^s}x[v_2v_1^{a-2}]x],
\]
equals some multiple of $[v_{a+p^s}xx]$, the actual coefficient being immaterial, and hence vanishes.
Consequently, we have
\[
0=[v_{p^s}[v_ayx]]=[v_{p^s}v_ayx]-[v_{p^s}xv_ay]-[v_{p^s}yv_ax]+[v_{p^s}xyv_a].
\]
Now we evaluate, in turn, each of the four summands.

By means of Lemma~\ref{lemma:technical} and relations~\eqref{eq:va-1} and~\eqref{eq:va} we obtain
\begin{align*}
[v_{p^s}v_a]
&=[v_{p^s}[v_2v_1^{a-2}]]
=\binom{a-2}{a-1-p^s}
[v_{a-1}v_2v_1^{p^s-1}]
-\binom{a-2}{a-p^s}
[v_{a}v_2v_1^{p^s-2}]\\
&=[v_{a-1}v_2v_1^{p^s-1}]
+[v_{a}v_2v_1^{p^s-2}]
=-2\,v_{a+p^s}+v_{a+p^s}=-v_{a+p^s},
\end{align*}
and hence $[v_{p^s}v_ayx]=-[v_{a+p^s}yx]$.
Similarly, we have
\begin{align*}
[v_{p^s}xv_ay]
=[v_{a-1}xv_2v_1^{p^s-1}y]
+[v_{a}xv_2v_1^{p^s-2}y]
=-[v_{a+p^s}xy].
\end{align*}
The third summand is
\[
[v_{p^s}yv_ax]=[v_{a-1}yv_2v_1^{p^s-1}x]
+[v_{a}yv_2v_1^{p^s-2}x]
=-[v_{a+p^s}yx],
\]
where we have used the additional relations
\[
[v_{a-1}yv_2]=-[v_{a+1}y],\qquad[v_ayv_2]=0.
\]
Finally, the fourth summand is
\[
[v_{p^s}xyv_a]=[v_{a-1}xyv_2v_1^{p^s-1}]
+[v_{a}xyv_2v_1^{p^s-2}]=-2[v_{a+p^s}xy].
\]
where we have used the further relations
\[
[v_{a-1}xyv_2]=-2[v_{a+1}xy],\qquad
[v_axyv_2]=0.
\]
It follows that
\[
0=[v_{p^s}[v_ayx]]=-[v_{a+p^s}xy],
\]
which concludes the proof.
\end{proof}

We can extract from the proof of Lemma~\ref{lemma:a-cong-1}
some information also on the case where $a=p^s+1$, namely that
\[
[v_{a+k}xx]=0 \quad\text{for}\quad 2 \leq k < p^s,
\]
but not for $k=p^s$.
Hence in this case $T$, if large enough, has further diamonds of infinite type in all degrees
$k(q-1)+1$ for $p^s+1<k<2p^s+1$.
This is one step in the direction of a proof of Theorem~\ref{thm:presentation}, for which we refer to~\cite{Avi:thesis}.

Some of the arguments and computations used in the course of the proof of Theorem~\ref{thm:L-dies}
have a wider applicability, as they only exploit the
structure of the thin algebra in a neighborhood of a diamond, besides the initial structure up to the second diamond.
This is clearly so for the technical Lemma~\ref{lemma:technical}, but also, for example, for
an argument preceding Lemma~\ref{lemma:a-notcong-1}.
That argument shows that if $T$ is a $(-1)$-algebra with second diamond in degree $2q-1$,
with $T_{m-q+1}$ a diamond of infinite type,
$T_m$ a diamond of finite type, and all components $T_{m-q+2}$ up to $T_{m-2}$ having
the same two-step centralizer as $T_2$, then $T_{m+q-1}$ is a diamond of infinite type,
provided $\dim(T)$ is large enough.

%----------------------------------------------------------------------------
\bibliography{References}

\end{document}